\documentclass[10pt]{amsart}
\usepackage{amssymb,amsmath, amsthm, amsfonts}

\newcommand{\rd}{\mathbb{R}^d}

\newcommand{\R}{\mathbb{R}}

\newcommand{\zed}{\mathbb{Z}}
\newcommand{\pd}{\frac{d(d+1)}{d^2-d+2}}

\newcommand{\pq}{L^p\to L^q}
\newcommand{\pqd}{L^p(\mathbb{R}^d)\to L^q(\mathbb{R}^d)}

\newcommand{\ct}{\mathcal{T}}

\newcommand{\tef}{\langle X\chi_E,\chi_F\rangle}

\newcommand{\tfe}{\langle \chi_E,X^*\chi_F\rangle}
\newcommand{\epo}{\epsilon_1}
\newcommand{\ept}{\epsilon_2}
\newcommand{\eto}{\eta_1}
\newcommand{\ett}{\eta_2}

\newtheorem{thm}{Theorem}

\newtheorem{lem}{Lemma}

\begin{document}
\title[Restricted X-ray transforms]{A note on restricted X-ray transforms}
\author{Norberto Laghi}
\address{School of Mathematics and Maxwell Institute of Mathematical
  Sciences, The University of Edinburgh, JCM
 Building, The King's Buildings, Edinburgh EH9 3JZ, Scotland}
\email{N.Laghi@ed.ac.uk}
\thanks{The author is supported by an EPSRC grant.}
\begin{abstract}We show how the techniques introduced in \cite{1} and \cite{2}
can be employed to derive endpoint $\pq$ bounds for the X-ray
transform associated to the line complex generated by the curve
$t\to (t,t^2,\ldots,t^{d-1}).$ Almost-sharp Lorentz space estimates
are produced as well.
\end{abstract}
\maketitle
\section{Introduction}
The purpose of this note is to give yet another application of the
far-reaching techniques of Christ (see \cite{1}), and in particular
to utilise the refinement provided in \cite{2} to establish strong
type $(p,q)$ bounds for the X-ray transform first studied
in \cite{4} and \cite{3}. Work on averages along curves using these
techniques  is currently being undertaken by Stovall (\cite{8});
however, the simple behaviour of the X-ray transform makes  
it a very natural object of study.
Thus, we shall be concerned with the operator
\begin{equation}\label{1}Xf(x)=
\int_I f(s,x_2+sx_1,x_3+sx_1^2,\ldots,x_d+sx_1^{d-1})ds,\end{equation}
 where $I\subset\R$ 
is a closed interval.\\
To characterise the set of $(p,q)$ such that $X:\pq,$ we let
$\Delta_d\subset [0,1]^2$ 
be the convex hull of the points $(1,1),(0,0)$ and 
$(p_d^{-1},q_d^{-1})$, where $p_d=\frac{d(d+1)}{d^2-d+2},\ 
q_d=\frac{d+1}{d-1}.$ We have the following.
\begin{thm} \[X:\pqd\iff (p^{-1},q^{-1})\in\Delta_d.\]
\end{thm}
The necessity of this condition can be found in \cite{6}.
We remark that this result was established in \cite{5} for $d=3,$
whilst \cite{6} and \cite{7} provide estimates in dimension $d=4$,
\cite{4} contains estimates in dimension $d=4,5$ and
\cite{3} contains estimates for general $d;$ further mixed-norm
estimates for the operator (\ref{1}) have been obtained as well, and 
we refer the reader to
\cite{3} and \cite{4} and the references therein for these results.
However 
the work of the aforementioned authors establishes only restricted 
weak-type bounds
at $(p_d,q_d)$ when $d\geq 4.$ We are able to establish the following
estimate, from which Theorem 1 follows by simple interpolation.
\begin{thm}For every $\epsilon>0,$ 
\[X:L^{p_d}(\rd)\to L^{q_d,p_d+\epsilon}.\]
\end{thm}
Here $L^{s,r}(\rd)$ denote the familiar Lorentz spaces;
since $p_d<q_d,$ it is clear that the strong $(p_d,q_d)$ bound for $X$
follows from Theorem 2.\\
The following example shows that $X$ maps $L^{p_d}$ to $L^{q_d,r}$
only if $r\geq p_d.$  Hence, with the possible exception of
$\epsilon,$ Theorem 2 is sharp in the scale of Lorentz
spaces.\\
For simplicity we shall suppose that 
$0\in I$ in the definition of (\ref{1}), otherwise we may translate 
things accordingly.
Let $\chi=\chi_{[-1,1]^d}$ be the
characteristic function of the cube centred at the origin of
sidelength $2.$ Now consider the nonisotropic dilations given by
\[\delta\circ y=(\delta y_1,\ldots,\delta^d y_d),\]
and let $\chi_k=\chi_k(y)=(k^{-1}\circ y).$ We define the function
\[f(x)=\sum_{k\geq N}\chi_k(x_1,x_2-k^2,\ldots,x_d-k^d),\quad N\gg
1,\] where the single elements in the sum have clearly disjoint
supports. Thus, if $A_k=supp(\chi_k(\cdot-(0,k^2,\ldots,k^d))$
\begin{multline*}\|f\|_{L^{p_d}}=\left(\sum_{k\geq
    N}|A_k|\right)^{\frac{d^2-d+2}{d(d+1)}}=
 \left(\sum_{k\geq
    N}k^{-d(d+1)/2}\right)^{\frac{d^2-d+2}{d(d+1)}}\approx\\
N^{(-d(d+1)/2+1)\frac{d^2-d+2}{d(d+1)}}=N^{(-1/2+1/{d(d+1)})(d^2-d+2)}.
\end{multline*}
However
\begin{multline*}Xf(x)=\sum_{k\geq N}\int_I\chi_k 
(s,x_2+sx_1-k,x_3+sx_1^2-k^2,\ldots,x_d+sx_1^{d-1}-k^d)ds
\geq\\ \sum_{k\geq N} \psi_{2k}(x_1)
\int_{s\ll k^{-1}} 
\chi_k 
(s,x_2+sx_1-k,x_3+sx_1^2-k^2,\ldots,x_d+sx_1^{d-1}-k^d)ds\gtrsim \\
\sum_{k\geq N}k^{-1}\psi_{2k}(x_1)
\psi_{{2k}^2}(x_2-k^2)\ldots
\psi_{{2k}^d}(x_d-k^d),\end{multline*} where
$\psi_j(t)=\chi_{[-1/j,1/j]}(t).$ Note that, again, all the functions
in the sum are characteristic functions of disjoint sets. Thus,
\begin{multline*}\|Xf\|_{L^{q_d,r}}\gtrsim \left[\sum_{k\geq N}
\left(k^{-1}|B_k|^{\frac{d-1}{d+1}}\right)^r\right]^{1/r}=
\left[\sum_{k\geq N}
\left(k^{-1}k^{-\frac{d(d+1)}{2}\frac{d-1}{d+1}}\right)^r\right]^{1/r}=\\
\left[\sum_{k\geq N}k^{-\frac{d^2-d+2}{2}r}\right]^{1/r}\approx
N^{(-\frac{d^2-d+2}{2}+1/r)}.
\end{multline*}
Hence, in order for boundedness to hold, we must have 
\begin{gather*}
N^{(-\frac{d^2-d+2}{2}+1/r)}\lesssim N^{(-1/2+1/{d(d+1)})(d^2-d+2)}
\Longrightarrow\\
-\frac{d^2-d+2}{2}+\frac1{r}\leq  -\frac{d^2-d+2}{2}+
\frac{d^2-d+2}{d(d+1)}\iff r\geq \pd\ ,\end{gather*}
as $N\to\infty.$\\
\emph{Notation.} Whenever we write (or we have written) $A\lesssim B$
for any two nonnegative quantities $A$ and $B,$ we shall mean that
there exists a striclty positive constant $c$ such that $A\leq cB;$
such constant is subject to change from line to line and even from
step to step.
 
\section{Preliminary statements}
We first summarise some of the results contained in \cite{3} that are
necessary to our arguments. Let $E\subset \rd_1,F\subset\rd_2$ 
be any two measurable
sets\footnote{We use the notation $\rd_i,\ i=1,2$ to stress the fact
that the sets $E$ and $F$ lie in different ambient spaces, albeit of
equal dimension $d.$}
and let
\begin{equation}\ct(E,F)=\tef=\tfe,\end{equation} where the dual
operator $X^*$ is given by
\[X^*g(x)=\int_{\R}g(t,x_2-x_1t,
x_3-x_1t^2,\ldots,x_d-x_1t^{d-1})dt.\]
The restricted weak-type $(p_d,q_d)$ estimate for (\ref{1}) then amounts
to prove
\[\tef\lesssim |E|^{1/{p_d}}|F|^{1/{{q_d}'}}.\] If one lets
\[\alpha=\ct(E,F)/{|F|},\quad \beta=\ct(E,F)/{|E|},\] then it suffices to
show that either
\begin{equation}\label{3}
|E|\gtrsim \alpha^d\beta^{d(d-1)/2}\quad\text{or}\quad 
|F|\gtrsim\alpha^{d-1}\beta^{(d^2-d+2)/2}\end{equation}
\begin{thm}(Christ-Erdogan, \cite{3}) Estimates (\ref{3}) hold.
\end{thm}
It is important to keep in mind the manner by which estimates
(\ref{3}) were established.
Define, for fixed $x,$ the maps 
\begin{align*}&\gamma(x,s)=(s,x_2+sx_1,x_3+sx_1^2,\ldots,x_d+sx_1^{d-1}),\\
&\gamma^*(x,t)=(t,x_2-x_1t,
x_3-x_1t^2,\ldots,x_d-x_1t^{d-1}),\end{align*} and the maps
$\Phi_j\equiv\Phi_{j,x}:\R^j\to \rd_1$ if $j$ is even, 
$\Phi_j\equiv\Phi_{j,x}:\R^j\to\rd_2$
if $j$ is odd, by letting
$\Phi_1(t_1)=\gamma^*(x,t_1),\Phi_2(t_1,s_1)=\gamma(\Phi_1(t_1),s_1)$
and further
\begin{align*}\Phi_{2k+1}(t_1,s_1,t_2,s_2,\ldots,t_{k+1})=&
\gamma^*(\Phi_{2k}(t_1,s_1,t_2,s_2,\ldots,t_k,s_k),t_{k+1}),\\
\Phi_{2k+2}(t_1,s_1,t_2,s_2,\ldots,t_{k+1},s_{k+1})=&
\gamma^(\Phi_{2k}(t_1,s_1,t_2,s_2,\ldots,s_k,t_{k+1},s_{k+1})\end{align*}
Further, we define maps $\Psi_j\equiv\Psi_{j,x}$
 by setting $\Psi_1(s_1)=\gamma(x,s_1),
\Psi_2(s_1,t_2)=\gamma^*(\Psi_1(s_1),t_2)$ and
\begin{align*}\Psi_{2k+1}(s_1,t_2,\ldots,s_k,t_{k+1},s_{k+1})=&
\gamma(\Psi_{2k}(s_1,t_2,\ldots,s_k,t_{k+1}),s_{k+1})\\
\Psi_{2k+2}(s_1,t_2,\ldots,t_{k+1},s_{k+1},t_{k+2})=&
\gamma^*(\Psi_{2k+1}(s_1,t_2,\ldots,s_k,t_{k+1},s_{k+1}),t_{k+2}),
\end{align*}
where $\Psi_j:\R^j\to \rd_2$ if $j$ is even, and $\Psi_j:\R^j\to\rd_1$
if $j$ is odd. 
\begin{lem}Consider the maps $\Phi_d$ and $\Psi_d$ and let
$J_{\Phi}$ and $J_{\Psi}$ be the determinants of the associated
Jacobian matrices. Then
\begin{align*}& J_{\Phi}=c_d\prod_{j=1}^k(s_j-s_{j-1})\prod_{1\leq
    j<\ell\leq k}(t_j-t_{\ell})^4\quad\text{for }d=2k,\\
& J_{\Phi}=c_d\prod_{j=1}^k(s_j-s_{j-1})\prod_{1\leq
    j<\ell\leq k}(t_j-t_l)^4\prod_{j=1}^k(t_j-t_{k+1})^2
\quad\text{for }d=2k+1,\end{align*} where we have set $x_1=s_0$ for
notational convenience, and
\[J_{\Psi}=c_d(t_{k+1}-t_1)\prod_{j=1}^{k-1}(s_{j+1}-s_j)\prod_{2\leq j<\ell\leq
 k}(t_j-t_{\ell})^4\prod_{j=2}^k(t_j-t_{k+1})^2\prod_{j=2}^k(t_j-t_1)^2\]
when $d=2k,$ whilst
\[ J_{\Psi}=c_d\prod_{j=1}^{k}(s_{j+1}-s_j)\prod_{2\leq j<\ell\leq
 k+1}(t_j-t_{\ell})^4\prod_{j=2}^{k+1}(t_j-t_1)^2\] when $d=2k+1.$ We
 set $x_1=t_1$ in the formulas characterising $J_{\Psi}.$
\end{lem}
\begin{proof} The formulas concerning $J_{\Phi}$ have been proven in
 \cite{3}. To compute $J_{\Psi}$ first observe that for even $d$
\begin{align*}&\Psi_d(s_1,t_2,\ldots,s_k,t_{k+1})=\\
&(t_{k+1},x_2+\sum_{j=1}^k(t_{j}-t_{j+1})s_j,x_3+\sum_{j=1}^k
(t_{j}^2-t_{j+1}^2)s_j,\ldots,
x_d+\sum_{j=1}^k(t_{j}^{d-1}-t_{j+1}^{d-1})s_j)\end{align*}
and for odd $d$
\begin{align*}&\Psi_d(s_1,t_2,\ldots,s_k,t_{k+1},s_{k+1})=\\
&(s_{k+1},x_2+\sum_{j=1}^k(t_{j}-t_{j+1})s_j+s_{k+1}t_{k+1},\ldots,
x_d+\sum_{j=1}^k(t_{j}^{d-1}-t_{j+1}^{d-1})s_j+s_{k+1}t_{k+1}^{d-1}).
\end{align*}
The corresponding Jacobian matrices  
${\partial \Psi_d}/{\partial(s,t)}$\footnote{It is important to keep
  in mind that $t_1$ is just a dummy variable.}
are then given by
\[
\left(\begin{matrix}1&-s_k& -2t_{k+1}s_k & \ldots&
  -(d-1)t_{k+1}^{d-2}s_k\\
0& t_k-t_{k+1}&t_k^2-t_{k+1}^2&\ldots & t_k^{d-1}-t_{k+1}^{d-1}\\
0 & s_k-s_{k-1}& 2t_k(s_k-s_{k-1})&
\ldots&(d-1)t_k^{d-2}(s_k-s_{k-1})\\
\vdots &\vdots&\vdots&\ddots&\vdots\\
0& t_2-t_3& t_2^2-t_3^2 &\ldots& t_2^{d-1}-t_3^{d-1}\\
0& s_2-s_1&2t_2(s_2-s_1)&\ldots&(d-1)t_2^{d-1}(s_2-s_1)\\
0& t_1-t_2& t_1^2-t_2^2&\ldots&
t_1^{d-1}-t_2^{d-1}\end{matrix}\right)\] for even $d,$ and by 
\[
\left(\begin{matrix}1&t_{k+1}& t_{k+1}^2 & \ldots&
  t_{k+1}^{d-1}\\
0& s_{k+1}-s_k&2t_{k+1}(s_{k+1}-s_k)&\ldots & (d-1)t_{k+1}^{d-2}(s_{k+1}-s_k)\\
0 & t_k-t_{k+1}& t_k^2-t_{k+1}^2&
\ldots&t_k^{d-1}-t_{k+1}^{d-1})\\
\vdots &\vdots&\vdots&\ddots&\vdots\\
0& t_2-t_3& t_2^2-t_3^2 &\ldots& t_2^{d-1}-t_3^{d-1}\\
0& s_2-s_1&2t_2(s_2-s_1)&\ldots&(d-1)t_2^{d-1}(s_2-s_1)\\
0& t_1-t_2& t_1^2-t_2^2&\ldots&
t_1^{d-1}-t_2^{d-1}\end{matrix}\right)\]
for odd $d.$\\
Let us first examine the case where $d=2k$ is even. Observe, just like in
\cite{3},
 that $J_{\Psi}$ must be a polynomial of degree $d(d-1)/2,$ and that
 the polynomial considered in the statement of the lemma has the same degree.
 Now, it suffices to factor out of the determinant all the terms of the form
 $(s_{j+1}-s_j),$ $j=1,\ldots, k-1,$ to obtain an expression involving only
 the $t$ variables. One can then prove that the determinant is divisible
 by the quadratic and quartic terms just like in \cite{3}, whilst the presence
 of the linear term can be seen by observing that by adding every other row
 starting from the bottom one, we obtain a matrix 
with  a row of the form
\[\left(\begin{matrix}t_1-t_{k+1}& t_1^2-t_{k+1}^2&\ldots&
t_1^{d-1}- t_{k+1}^{d-1}\end{matrix}\right),\]
where all entries have the common factor $(t_1-t_{k+1}).$ 
One may then prove that the constant
$c_d\neq 0,$ again as in \cite{3}.\\
The case of odd $d=2k+1$ is even simpler; again  we may factor all
the terms of the form $(s_{j+1}-s_j),$ $j=1,\ldots,k$ to obtain an expression involving
only the $t$ variables. However, this is completely analogous to the case treated in
\cite{3}, and the same can be said about the constant $c_d.$ 
\end{proof}
We conclude this section with a lemma that can be seen as the natural analogue
of Lemma 8.1 in \cite{2}.
\begin{lem}Let $E,E'\subset\rd_1,\ G\subset\rd_2$ be measurable sets of finite measure, and suppose
that $X\chi_{E'}(x)\geq\delta_1$ for all $x\in G.$ Then
\begin{equation}\label{4}|E'|\gtrsim \delta_1^2(\ct(E,G)/{|G|})^{d-2}(\ct(E,G)/{|E|})^{d(d-1)/2}.\end{equation}
Further, let $F,F'\subset\rd_2,\ H\subset\rd_1$ be measurable sets of finite measure, and suppose
$X^*\chi_{F'}(y)\geq\delta_2$ for all $y\in H.$ Then
\begin{equation}\label{5}|F'|\gtrsim \delta_2^d(\ct(H,F)/{|F|})^{d-1}(\ct(H,F)/{|F|})^{(d^2-d+2)/2-d}.\end{equation}
\end{lem}
\begin{proof}
We first prove (\ref{4}), by splitting the argument in the two cases of even $d$ and odd $d.$\footnote{This is why we shall need to utilise the formulae we derived for the maps $\Psi_d,$ as well as the formulae
for $\Phi_d.$} To simplify notation, we shall write $z=(z_1,\ldots,z_{m-1},z_m)=(\hat{z},z_{m})\in\R^m$ for any variable $z$ and appropriate $m\in\zed_+.$ \\
\textbf{Case} $d=2k.$  By using the method of refinements developed in \cite{1}, we may find
a point $x_0\in E$ and a sequence of sets $\Omega_j\subset \R^j,$ $j=1,\ldots,d$ satisfying
\begin{enumerate}\item for each $j,$ $\Omega_{j+1}\subset\Omega_j\times\R,$
\item $|\Omega_1|\gtrsim \ct(E,G)/{|E|},$
\item for even $j$, for each point $\omega\in\Omega_j,\ |\{t\in\R:(\omega,t)\in\Omega_{j+1}\}|\gtrsim
\ct(E,G)/{|E|},$
\item for odd $j\neq d-1$, for each point $\omega\in\Omega_j,\ |\{t\in\R:(\omega,t)\in\Omega_{j+1}\}|\gtrsim \ct(E,G)/{|G|},$
\item for $j=d-1,$ for each point $\omega\in\Omega _{d-1},\ |\{s\in\R:(\omega,s)\in\Omega_d\}|\gtrsim\delta_1,$
\item $\Phi_{j,x_0}(\Omega_j)\subset E$ for even $j,$ $\Phi_{j,x_0} (\Omega_j)\subset F$ for odd $j,$ and
$\Phi_{d,x_0}(\Omega_d)\subset E'.$ \end{enumerate}
Thus, by Bezout's theorem (see \cite{1},\cite{3}) we have the lower bound
\begin{multline*}|E'|\gtrsim\Phi_{d,x_0}(\Omega_d)\gtrsim\int_{\Omega_d}|J_{\Phi}(s,t)|dsdt=\\
\int_{\Omega_d}\prod_{j=1}^k|s_j-s_{j-1}|\prod_{1\leq
    j<\ell\leq k}|t_j-t_{\ell}|^4dsdt\gtrsim \\
  \delta_1^2\int_{\Omega_{d-1}}\prod_{j=1}^{k-1}|s_j-s_{j-1}|\prod_{1\leq
    j<\ell\leq k}|t_j-t_{\ell}|^4d\hat{s}dt\gtrsim\\ 
\delta_1^2(\ct(E,G)/{|G|})^{d-2}(\ct(E,G)/{|E|})^{d(d-1)/2}.
    \end{multline*}
\textbf{Case} $d=2k+1.$ Here the method of refinements gives us     
a point $y_0\in F$ and a sequence of sets $\Omega_j\subset \R^{j-1},$ $j=2,\ldots,d+1$ satisfying
\begin{enumerate}\item for each $j,$ $\Omega_{j+1}\subset\Omega_j\times\R,$
\item $|\Omega_2|\gtrsim \ct(E,G)/{|G|},$
\item for odd $j\neq d$, for each point $\omega\in\Omega_j,\ |\{s\in\R:(\omega,s)\in\Omega_{j+1}\}|\gtrsim
\ct(E,G)/{|G|},$
\item for even $j$, for each point $\omega\in\Omega_j,\ |\{t\in\R:(\omega,t)\in\Omega_{j+1}\}|\gtrsim \ct(E,G)/{|E|},$
\item for $j=d$ for each point $\omega\in\Omega _{d},\ |\{s\in\R:(\omega,s)\in\Omega_{d+1}\}|\gtrsim\delta_1,$
\item $\Psi_{j,y_0}(\Omega_{j+1})\subset E$ for odd $j,$ $\Psi_{j,y_0}(\Omega_{j+1})\subset F$ for
even $j,$ and
$\Psi_{d,y_0}(\Omega_{d+1})\subset E'.$ \end{enumerate} Again, by
Bezout's theorem,
\begin{multline*}|E'|\gtrsim\Psi_{d,y_0}(\Omega_{d+1})\gtrsim\int_{\Omega_{d+1}}|J_{\Phi}(s,t)|dsdt=\\
\int_{\Omega_{d+1}}\prod_{j=1}^{k}|s_{j+1}-s_j|\prod_{2\leq j<\ell\leq
 k+1}|t_j-t_{\ell}|^4\prod_{j=2}^{k+1}|t_j-t_1|^2dsdt\gtrsim\\
\delta_1^2
\int_{\Omega_{d}}\prod_{j=1}^{k-1}|s_{j+1}-s_j|\prod_{2\leq j<\ell\leq
 k+1}|t_j-t_{\ell}|^4\prod_{j=2}^{k+1}|t_j-t_1|^2d\hat{s}dt\gtrsim\\
 \delta_1^2(\ct(E,G)/{|G|})^{d-2}(\ct(E,G)/{|E|})^{d(d-1)/2}.\end{multline*}
We now turn to the proof of  (\ref{5}).\\
\textbf{Case}  $d=2k+1.$ Here we may apply the previous method as in the case of even $d$
for (\ref{4}); however, now properties (5) and (6) should now be\\
(5) for $j=d-1,$ for each point $\omega\in\Omega _{d-1},\ |\{t\in\R:(\omega,t)\in\Omega_d\}|\gtrsim\delta_2,$\\
(6) $\Phi_{j,x_0}(\Omega_j)\subset E$ for even $j,$ $\Phi_{j,x_0} (\Omega_j)\subset F$ for odd $j,$ and
$\Phi_{d,x_0}(\Omega_d)\subset F'.$\\
The lower bound one gets thanks to Bezout's theorem is now
\begin{multline*}|F'|\gtrsim |\Phi_{d,x_0}(\Omega_d)|\gtrsim \int_{\Omega_d} |J_{\Phi}(s,t)|dsdt=\\
\int_{\Omega_d} \prod_{j=1}^k|s_j-s_{j-1}|\prod_{1\leq
    j<\ell\leq k}|t_j-t_l|^4\prod_{j=1}^k|t_j-t_{k+1}|^2dsdt\gtrsim\\
\delta_2^{2k+1}\int_{\Omega_{d-1}} \prod_{j=1}^k|s_j-s_{j-1}|\prod_{1\leq
    j<\ell\leq k}|t_j-t_l|^4dsd\hat{t} \gtrsim\\ 
\delta_2^d (\ct(H,F)/{|F|})^{d-1}(\ct(H,F)/{|H|})^{(d^2-d+2)/2-d}.
\end{multline*}
\textbf{Case}  $d=2k.$ 
Here we may apply the method of refinements as in the case of odd $d$
for (\ref{4}); now conditions (5) and (6) are\\
(5) for $j=d$ for each point $\omega\in\Omega _{d},\ |\{t\in\R:(\omega,t)\in\Omega_{d+1}\}|\gtrsim\delta_2,$\\
(6) $\Psi_{j,y_0}(\Omega_{j+1})\subset E$ for odd $j,$ $\Psi_{j,y_0}(\Omega_{j+1})\subset F$ for
even $j,$ and
$\Psi_{d,y_0}(\Omega_{d+1})\subset F'.$ Thus, using Bezout's theorem
once more, we have
\[|F'|\gtrsim|\Psi_{d,y_0}(\Omega_{d+1}|\gtrsim\int_{\Omega_{d+1}}|J_{\Psi}(s,t)|dsdt
\gtrsim\]\begin{align*}
\int_{\Omega_{d+1}}|t_{k+1}-t_1|\prod_{j=1}^{k-1}|s_{j+1}-s_j|&
\prod_{2\leq j<\ell\leq
 k}|t_j-t_{\ell}|^4\prod_{j=2}^k|t_j-t_{k+1}|^2\prod_{j=2}^k|t_j-t_1|^2
dsdt    \gtrsim\\
 \delta_2^{2k}\int_{\Omega_d}\prod_{j=1}^{k-1}|s_{j+1}-s_j&|\prod_{2\leq j<\ell\leq
 k}|t_j-t_{\ell}|^4\prod_{j=2}^k|t_j-t_{k+1}|^2dsd\hat{t}\gtrsim\\ 
&\delta_2^d (\ct(H,F)/{|F|})^{d-1}(\ct(H,F)/{|H|})^{(d^2-d+2)/2-d}.   
 \end{align*}
 \section{Strong type estimates}
The purpose of this section is to show how the arguments in \cite{2}
can be utilised to obtain the statement of Theorem 2; naturally,
we shall have to make suitable modifications, the main one being the use
of Lemma 2. We are aiming to show that
\begin{equation}\label{6}\left|\langle Xf,g\rangle\right|\lesssim
  \|f\|_{L^{p_d}}
\|g\|_{L^{{q_d}',r'}},\quad r'<{p_d}',\end{equation}
 which naturally implies $X:L^{p_d}(\rd)\to L^{q_d,r}(\rd)$ for
$r>p_d.$ For the sake of notational simplicity, from now on we shall relabel
$p\equiv p_d,\ q\equiv q_d,$ as we shall only deal with inequality
(\ref{6}) in this section.\\
As pointed out in \cite{2}, it suffices to consider $f,g$ of the form
$f=\sum_{k\in\zed}2^k\chi_{E_k},\ g=\sum_{j\in\zed}2^j\chi_{F_j}$
were the sets $E_k$'s are pairwise disjoint and so are the $F_j$'s;
the indices $k,j$ are completely independent of each other. The key
step is to show that
\begin{equation}\label{7}\left|\langle Xf,g\rangle\right|\lesssim \|f\|_p\|g\|_{q'}\quad\text{if }
f=\sum_{k\in\zed}2^k\chi_{E_k}\ \text{ and }\
g\equiv\chi_F\quad\text{for a single set }F,\end{equation}
and its counterpart\footnote{This is utterly redundant when the
operator in question is (essentially) self-adjoint, but the X-ray
transform does not have this property.}
\begin{equation}\label{8}\left|\langle Xf,g\rangle\right|\lesssim \|f\|_p\|g\|_{q'}\quad\text{if }
f\equiv\chi_E\quad\text{for a single set }E,\quad\text{and
}g=\sum_{j\in\zed}2^j\chi_{F_j}.\end{equation}
We follow the scheme of \cite{2} to prove (\ref{7}). Let
$\epo,\eto\in(0,1/2]$ be arbitrary and normalise the $p$ norm of $f$
  by setting $\sum_k 2^{kp} |E_k|=1.$ Suppose
\[|E_k|\approx\eto 2^{-kp}\quad\text{for all }k,\quad \ct(E_k F)\approx
\epo |E_k|^{1/p}|F|^{1/q'}\quad\text{for all }k.\]
Then the number $M_1$ of indices $k$ is finite and $M_1\eto\lesssim
1.$ Further, assume that any two indices $k_1,k_2$ in the sum satisfy
$|k_1-k_2|\geq A\log(1/{\epo}),$\footnote{This is done by simply
splitting the sum into $O(A\log{(1/{\epo})})$ sums; the logarithmic
factor that is lost will not affect the estimates in a crucial way.}
and define the sets
\[G_k=\left\{x\in F:X{\chi_{E_k}}(x)\geq c_0\epo
|E_k|^{1/p}|F|^{1/{q'}-1}\right\},\] where the constant $c_0>0$ is
 chosen sufficiently small to have $\ct(E_k,F\setminus G_k)\leq \frac12
 \ct(E_k F),$ so that $\mathcal{T}(E_k,G_k)\approx\ct(E_k,F).$
Since $\ct(E_k,G_k)\lesssim |E_k|^{1/p}|G_k|^{1/{q'}},$ this implies
\begin{equation}\label{9} |G_k|\gtrsim \epo^{q'}|F|.\end{equation}
A simple observation\footnote{So far we have only described the
arguments in \cite{2}, which we have included for the sake of completeness;
all the details can be found in that paper.} then shows that one has
the dichotomy
\begin{equation}\label{10}\text{either }\sum_{k\in\zed}|G_k|\lesssim
|F|,\quad\text{or}\end{equation} there exists indices $k_1,k_2,$ $k_1\neq k_2$
so that
\begin{equation}\label{11} |G_{k_1}\cap G_{k_2}|\gtrsim\epo^{2q'}|F|.
\end{equation} We first wish to show that (\ref{11}) cannot hold; we
shall then complete the proof as in \cite{2}.\\
Arguing by contradiction, assume that (\ref{11}) does hold; we 
start by applying (\ref{4}) of Lemma 2 with $E=E_{k_1}, E'=E_{k_2},
G=G_{k_1}\cap G_{k_2}$ and $\delta_1\approx \epo |E'|^{1/p}|F|^{-1/q}.$
We also have that $X{\chi_E}\geq c_0\epo|E|^{1/p}|F|^{-1/q}$ at each
point of $G$ and thus
\[\ct(E,G)\gtrsim\epo|E|^{1/p}|F|^{-1/q}|G|.\] By Lemma 2 one may
conclude
\begin{multline*}|E'|\gtrsim (\epo |E'|^{1/p}|F|^{-1/q})^2
(\epo |E|^{1/p}|F|^{-1/q})^{d-2}\\(\epo
  |E|^{1/p-1}|F|^{-1/q}|G|)^{d(d-1)/2}
\gtrsim\\ \epo^{d+q'd(d-1)/2}|E'|^{2/p}|E|^{(d-2)/p}|F|^{-d/q}
|F|^{d(d-1)/{2q'}}=\\ 
\epo^{d+q'd(d-1)/2}|E'|^{2/p}|E|^{(d-2)/p},\end{multline*}
where we have used (\ref{11}) and the actual expressions for $(p,q).$
After a bit of algebra one then reaches the conclusion
\[|E'|\lesssim\epo^{-\varphi}|E|,\quad\text{for some }\varphi>0.\]
From here, since $|E|=|E_{k_1}|\approx\eto 2^{-k_1p}$ and 
$|E'|=|E_{k_2}|\approx\eto 2^{-k_2 p}$ and the fact that the roles of
$E$ and $E'$ can be interchanged, one obtains that $|k_1-k_2|\lesssim 
\log(1/{\epo}),$ a contradicion to the assumption 
$|k_1-k_2|\geq A\log(1/{\epo})$ if $A$ is chosen sufficiently large.
\vspace{1mm}\\
Hence (\ref{10}) holds and we may now conclude the argument.  
We have
\begin{multline*}\sum_{k\in\zed}2^k\ct(E_k,F)\approx\sum_{k\in\zed}2^k
\ct(E_k,G_k)
\lesssim\\ (\sum_{k\in\zed}2^{kq}|E_k|^{q/p})^{1/q}
(\sum_{k\in\zed}|G_k|)^{1/{q'}}\lesssim \\
(\sum_{k\in\zed}2^{kp}|E_k|2^{k(q-p)}|E_k|^{q/p-1})^{1/q}|F|^{1/{q'}}\leq\\
\max_k
(2^{kp}|E_k|)^{(1/p-1/q)}|F|^{1/{q'}}\lesssim\eta^{1/p-1/q}|F|^{1/{q'}},
\end{multline*}
where $1/p-1/q>0$ and we used $\sum_{k\in\zed}2^{kp}|E_k|=1.$
However, since the number of
indices $k$ in the sum is $M_1\lesssim\eto^{-1},$ one may also argue
that
\begin{multline*}\sum_{k\in\zed}2^k\ct(E_k,F)\approx\sum_{k\in\zed}2^k\eto
|E_k|^{1/p}|F|^{1/{q'}}\lesssim\\ \epo M_1\eto^{1/p}|F|^{1/{q'}}=\epo
\eto^{-1/{p'}}|F|^{1/{q'}}.\end{multline*}
If we now recall the assumption $|k_1-k_2|\geq A\log(1/{\epo})$  and
  retain the normalisations in $\epo,\eto$ we have
$\langle
  Tf,\chi_F\rangle\lesssim\log(1/{\epo})\min(\eto^{1/p-1/q},\epo
\eto^{-1/{p'}})|F|^{1/{q'}}.$ Thus,
\begin{equation}\label{12}\langle
  Tf,\chi_F\rangle\lesssim\min
  (\epo^a,\eto^b)\|f\|_p|F|^{1/{q'}}\end{equation}
for positive $a,b$ and all $f,F$ subject to the normalisations in
  $\epo,\eto.$ Summing over dyadic values of $\eto$ we have
\begin{equation}\label{13}\langle
  Tf,\chi_F\rangle\lesssim\epo^a\|f\|_p|F|^{1/{q'}}\end{equation}
where now $f,F$ are only subject to $\epo$ normalisations. Summing
again over dyadic values of $\epo$ gives (\ref{7}), although it is 
equation (\ref{13}) we shall use to prove the strong type bounds.
\vspace{1mm}\\
We now give an outline of the argument needed to prove (\ref{8}).
Again, let $\ept,\ett\in(0,1/2]$ be arbitrary, and normalise the 
$q'$ norm of $g$ by setting $\sum_{j\in\zed}2^{jq'}|F_j|=1.$
Suppose\[|F_j|\approx\ett 2^{-kq'}\quad\text{for all
}j,\quad\ct(E,F_j)\approx\ept |E|^{1/p}|F_j|^{1/q'}\quad\text{for all }j.\]
We define $M_2$ as the number of indices $j$ in the sum, and again
$M_2\ett\lesssim 1;$ further, we shall split the sum in
$O(\log(1/{\ett}))$ sums. If we define
\[H_j=\left\{x\in E:X^*\chi_{F_j}
\gtrsim d_0\ept|F_j|^{1/{q'}}|E|^{1/p-1}\right\}\]
where $d_0$ is to be chosen sufficiently small so that
\[\ct(H_j,F_j)\approx \ct(E,F_j).\] Proceeding as in the proof of
(\ref{7}) one deduces that $|H_j|\gtrsim\ept^p |E|,$
and the new dichotomy becomes that either
$\sum_{j\in\zed}|H_j|\lesssim|E|,$ 
or there exist $j_1,j_2$ with $j_1\neq j_2$ so that 
$|H_{j_1}\cap H_{j_2}|\gtrsim\ept^{2p}|E|.$ Again, the key step is now
to show that the latter can't happen by applying (\ref{5}) of Lemma 2
in the following manner; set $F=F_{j_1}, F'=F_{j_2},
H=H_{j_1}\cap H_{j_2},$ and $\delta_2\approx \ept |F'|^{1/{q'}}|E|^{-1/{p'}}.$
Further, we have $X^*\chi_{F}\gtrsim\ept |E|^{-1/{p'}}|F|^{1/{q'}}$ at
every point of $H,$ hence
\[\ct(H,F)\gtrsim \ept |E|^{-1/{p'}}|F|^{1/{q'}}|H|.\] By Lemma 2 we
can now conclude
\begin{multline*}|F'|\gtrsim (\ept |F'|^{1/{q'}}|E|^{-1/{p'}})^d
(\ett |F|^{1/{q'}-1}|E|^{-1/{p'}}|H|)^{d-1}\\
(\ett |F|^{1/{q'}}|E|^{-1/{p'}})^{(d^2-d+2)/2-d}\gtrsim\\
\ept^{\psi}|F'|^{d/{q'}}|F|^{-(d-1)/q+(d^2-d+2)/{2q'}},\end{multline*}
for some $\psi>0,$ where we used that $|H|\gtrsim\ept^{2p}|E|.$ The
same rearrengement as before then provides the desired contradiction
and shows that $\sum_{j\in\zed}|H_j|\lesssim |E|.$ Inequality
(\ref{8}) can then be proven just like inequality (\ref{7}).
\vspace{1mm}\\
\emph{Conclusion of the proof.} Now let 
$f=\sum_{k\in\zed}2^k\chi_{E_k},\ g=\sum_{j\in\zed}2^j\chi_{F_j},$
assume $\|f\|_p=\|g\|_{q'}=1,$ and let $\ept,\ett\in(0,1/2].$
We shall suppose $|F_j|\approx\ett 2^{-kq'}$ for all $j$ with
$|F_j|>0.$ Then we consider the sum $\sum_{j,k}^*\ct(E_k,F_j)$ where
the $*$ indicates that the sum is taken only with respect to $j,k$ or
pairs $(j,k)$ with $\ct(E_k,F_j)\approx\ept |E_k||F_j|^{1/{q'}}.$
Again, one assumes $|j_1-j_2|\geq B\log(1/{\ept}).$ The proof of
inequality (\ref{8}) gives us, for each pair $(j,k)$ sets
$H_{j,k}\subset E_k$ so that $\ct(E_k,F_j)\approx\ct(H_{j,k},F_j)$ and
$\sum_j^*|H_{j,k}|\lesssim |E_k|.$ Hence
\begin{multline*}\sum_{j,k}^*2^j2^k\ct(E_k,F_j)\lesssim\sum_{j,k}2^j2^k
\ct(H_{j,k},F_j)=\\
\sum_j 2^j\langle X(\sum_k^*
2^k\chi_{H_{j,k}}),\chi_{F_j}\rangle\lesssim
2^j |F_j|^{1/{q'}}(\sum_k^* 2^{kp}|H_{j,k}|)^{1/p},\end{multline*}
where in the last step one uses inequality (\ref{7}). By H\"{o}lder's
inequality this last quantity is controlled by
\begin{equation}\label{14}(\sum_j 2^{jp'}|F_j|^{p'/{q'}})^{1/p'}
(\sum_j\sum_k^*2^{kp}|H_{j,k}|)^{1/p}\lesssim
\ett^{1/{q'}-1/{p'}}(\sum_k^*2^{kp}|E_k|)^{1/p}.\end{equation}
On the other hand, one may use the alternative bound
\[\sum_{j,k}^*\ct(E_k,F_j)\lesssim \epo^a\sum_j 2^j|F_j|^{1/{q'}}
(\sum_k 2^{kp}|E_k|)^{1/p}\leq \epo^a M_2\ett^{1/{q'}}\lesssim\epo^a 
\ett^{-1/q},\] where in the first step inequality (\ref{13}) has been used.
Now, summing over dyadic values of $\epo$ and $\ett$
gives the strong $(p,q)$ bound; the Lorentz space bound may be
obtained by observing that the first term of (\ref{14}) may be
controlled
by $\sum_j 2^{jr'}|F_j|^{r'/{q'}}$ if $r'<p'.$ This implies that
$r>p,$ giving the conclusion ot Theorem 2. 
\section{Final remarks}
The material presented in this paper is an interesting application
of the techniques first introduced in \cite{1}, and then further
developed in \cite{2}. Whilst a number of results have been proven
by utilising these ideas, it is not yet clear to which extent these 
techniques can be applied, although it is perhaps fair to say that
the $\pq$ regularity of many interesting operators may be studied this
way. In \cite{2} Christ has already shown that one need not be restricted
to studying averages along curves, but may consider submanifolds of $\rd$
of higher dimension, specifically the paraboloid. Further, the fact
that strong type estimates may be established by exploiting the
Lorentz ``smoothing'' that these objects present at the endpoints 
suggests that endpoint estimates may be established as well, at least 
in the case of translation-invariant operators. The work of Stovall in
\cite{8}, as well as the simple application we gave in this article 
certainly raise hope that this may indeed be possible.
 
\end{proof}

\end{document}